\documentclass[a4paper, 10pt]{article}

\usepackage{amssymb}
\usepackage{amsmath}
\usepackage[dvips]{graphicx}% to include graphics
\usepackage[all]{xy}        % Xy-pic
\usepackage{alltt}          % print verbatim + LaTeX code
\usepackage{shortvrb}       % shorthand for \verb with \MakeShortVerb{} and \DeleteShortVerb{|}
\usepackage{amsthm}         % AMS theorem package
\usepackage{url}            % insert Internet addresses with \url
\usepackage{fancyhdr}       % fancy headers
\usepackage{multirow}       % use \multirow{nrows}{width}{text} command in tables/arrays
\newcommand{\EG}{\ensuremath{\underline{E}G}}
\newcommand{\E}[1]{\ensuremath{\underline{E}#1}}   %enable the command to be called in both math & text mode
\newcommand{\F}{\ensuremath{\mathfrak{F}}}

\newcommand{\orbitcat}{\ensuremath{\mathcal{O}_\mathfrak{F}G}}
\newcommand{\Fin}{\ensuremath{\mathfrak{Fin}}}

\newcommand{\stab}{\textup{stab}}
\newcommand{\ttimes}{\!\times\!}

\newcommand{\SL}{\ensuremath{SL(3,\mathbb{Z})}}
\newcommand{\GL}{\ensuremath{GL(3,\mathbb{Z})}}
\newcommand{\SLn}{\ensuremath{SL(n,\mathbb{Z})}}
\newcommand{\Z}{\ensuremath{\mathbb{Z}}}
\newcommand{\R}{\ensuremath{\mathbb{R}}}
\newcommand{\RR}{\ensuremath{\mathcal{R}}}

\numberwithin{equation}{section} %Change equation numbering into section.number

\newtheorem{proposition}{Proposition}
\newtheorem{theorem}[proposition]{Theorem}
\newtheorem{corollary}[proposition]{Corollary}

\theoremstyle{definition} \newtheorem{definition}{Definition}
\theoremstyle{remark} \newtheorem{remark}[proposition]{Remark}%

\title{Bredon homology and equivariant $K$-homology of $\SL$}
\author{Rub\'en S\'anchez-Garc\'\i{}a\thanks{Funded by the EPSRC and the School of
Mathematics, University of Southampton}}
\date{\today}

\begin{document}

\maketitle

\begin{abstract}
    \noindent We obtain the equivariant $K$-homology of the classifying
    space $\underline{E}\SL$ from the computation of its Bredon homology
    with respect to finite subgroups and coefficients in the representation
    ring. We also obtain the corresponding results for $\GL$. Our
    calculations give therefore the topological side of the
    Baum-Connes conjecture for these groups.
\end{abstract}

\section{Introduction}
Consider a discrete group $G$. The Baum-Connes conjecture
\cite{BCH94} identifies the $K$-theory of the reduced
$C^*$-algebra of $G$, $C^*_r(G)$, with the equivariant
$K$-homology of a certain classifying space associated to $G$.
This space is called the classifying space for proper actions,
written \EG{}. The conjecture states that a particular map between
these two objects, called the assembly map,
    $$
    \mu_i : K_i^G(\E G) \longrightarrow K_i(C_r^*(G)) \quad
    i \ge 0\, ,
    $$
is an isomorphism. Here the left hand side is the equivariant
$K$-homology of \EG{} and the right hand side is the $K$-theory of
$C^*_r(G)$. The conjecture can be stated more generally
\cite[Conjecture 3.15]{BCH94}.

The equivariant $K$-homology and the assembly map are usually
defined in terms of Kasparov's $KK$-theory. For a discrete group
$G$, however, there is a more topological description due to Davis
and L\"uck \cite{DL98}, and Joachim \cite{Joa03} in terms of
spectra over the orbit category of $G$. We will keep in mind the
topological picture of the Baum-Connes conjecture (see Mislin's
notes in \cite{MislinValette:ProperGroupActions}).

Part of the importance of this conjecture is due to the fact that
it is related to many other relevant conjectures in different
areas of mathematics \cite{MislinValette:ProperGroupActions}.
Nevertheless, the conjecture itself allows the computation of the
$K$-theory of $C^*_r(G)$ from the $K^G$-homology of \EG{}. In
turn, this $K$-homology can be achieved by means of the Bredon
homology of \EG, as we explain later.

The Baum-Connes conjecture has been proved for some large families
of groups, yet remains unsolved in general. In particular, Higson
and Kasparov \cite{HK97} proved the conjecture for groups having
the Haagerup property (or a-T-menable), that is, groups which
admits a metrically proper isometric action on some affine Hilbert
space. On the other side, a group has Kazhdan's property $T$ if
every isometric action of $G$ on an affine Hilbert space has a
fixed point. Thus, infinite groups with the property $T$ do not
have the Haagerup property. There are not many infinite discrete
groups with the property $T$ for which the Baum-Connes conjecture
has been proved (the first examples are due to Lafforgue
\cite{Laf98}).

Consequently, the group \SL{} becomes relevant in this context
since the Baum-Connes conjecture is unknown for $\SLn$, $n \ge 3$
and these groups have property $T$. On the other hand, the
Baum-Connes assembly map is known to be injective for \SLn{} (in
general, for all closed subgroups of a Lie group with a finite
number of connected components, see \cite[\S 7]{BCH94}). Finally,
note that there are counterexamples to the Baum-Connes conjecture
for groupoids that can be constructed from \SL{}, and more
generally for a discrete group with property $T$ and such that the
Baum-Connes map is injective (\cite[p.~338]{HLS02}).

In this paper we obtain the equivariant $K$-homology of \EG{} for
$G = \SL$ and $G = \GL$ from the computation of its Bredon
homology. The results amount to the topological side of the
Baum-Connes conjecture, $K_i^G(\EG)$. We state the results here.

\begin{theorem}\label{thm:theorem1} The Bredon homology of $\underline{E}\SL$ with coefficients in the
representation ring is
$$ H_i^\Fin\left(\underline{E}\SL; \mathcal{R}\right) = \left\{
\begin{array}{cl} \Z^{\oplus 8} & i = 0\\ 0 & i \neq 0 \end{array}
\right. \, . $$
\end{theorem}

\begin{corollary}\label{cor:corollary2} The equivariant $K$-homology of
$\underline{E}G$ for $G = \SL$ is
$$ K_0^G\left( \EG \right) = \Z^{\oplus 8}\, , \quad
K_1^G\left( \EG \right) = 0 \, .$$
\end{corollary}
The results for $\GL$ follow from a K\"unneth formula for Bredon
homology since $\GL = \SL \times \Z/2\Z$.

\begin{theorem}\label{thm:theorem3} The Bredon homology of $\underline{E}\GL$ with coefficients in the
representation ring is
$$ H_i^\Fin\left(\underline{E}\GL; \mathcal{R}\right) = \left\{
\begin{array}{cl} \Z^{\oplus 16} & i = 0\\ 0 & i \neq 0 \end{array}
\right. \, . $$
\end{theorem}

\begin{corollary}\label{cor:corollary4} The equivariant $K$-homology of
$\underline{E}G$ for $G = \GL$ is
$$ K_0^G\left( \EG \right) = \Z^{\oplus 16}\, , \quad
K_1^G\left( \EG \right) = 0 \, .$$
\end{corollary}

We start with a brief review on classifying spaces, Bredon
homology and equivariant $K$-homology; then we describe a model of
$\underline{E}\SL$; in the last section we compute its Bredon
homology and prove the theorems above. These results are part of
the author's PhD thesis \cite[Chapter 4]{SanchezGarciaPhD}.

I would like to thank my PhD supervisor Ian Leary for suggesting
this problem and support throughout the work. Thanks also to
Christine Vespa and Samuel W\"uthrich for useful comments on a
previous version of this paper.

\section{Preliminaries}\label{section:Preliminaries}
\subsection{Classifying space for proper actions}\label{section:ClassifSpaces}
Let $G$ be a discrete group. A $G$-CW-complex is a CW-complex with
a $G$-action permuting the cells and such that if a cell is sent
to itself, it is done by the identity map. We call the $G$-action
\emph{proper} if all cell stabilizers are finite subgroups of $G$.

\begin{definition}
A \emph{model for \EG{}} is a proper $G$-CW-complex $X$ such that
for any proper $G$-CW-complex $Y$ there is a unique $G$-map $Y
\rightarrow X$, up to $G$-homotopy equivalence.
\end{definition}
\noindent One can prove that a proper $G$-CW-complex $X$ is a
model of \EG{} if and only if the subcomplex of fixed points $X^H$
is contractible for each $H \le G$ finite. It can be shown that
classifying spaces proper actions always exists. They are clearly
unique up to $G$-homotopy equivalence. We write $\underline{B}G$
for the quotient $\EG/G$. Note that for \emph{free} actions
instead of proper, we recover the definition of $EG$, whose
quotient $BG$ is the classifying space for principal $G$-bundles.
See \cite[\S 2]{BCH94} or \cite{LN01} for details and more
information on classifying spaces.

\subsection{Bredon (co)homology}\label{section:BredonHomology}
Given a group $G$ and a family \F{} of subgroups, we will write
\orbitcat{} for the \emph{orbit category}. The objects are left
cosets $G/K$, $K \in \F$, and morphisms the $G$-maps $\phi: G/K
\rightarrow G/L$. Such a $G$-map is uniquely determined by its
image $\phi(K) = gL$, and we have $g^{-1}Kg \subset L$.
Conversely, such $g \in G$ defines a $G$-map.

A \emph{left} (resp.~\emph{right}) \emph{Bredon module} is a
covariant (resp.~contravariant) functor from $\orbitcat$ to the
category of abelian groups. Bredon modules form a category, which
is abelian, and we can use homological algebra to define Bredon
homology (see \cite[pp.~7-10]{MislinValette:ProperGroupActions}).
Nevertheless, we give now a practical definition.

Consider a $G$-CW-complex $X$, a family $\F$ of subgroups of $G$
containing all cell stabilizers, and a left Bredon module $M$. The
\emph{Bredon homology groups} $H^\F_i \left(X;M\right)$ are
obtained as the homology of the following chain complex $(C_*,
\partial_*)$. Let $\{e_\alpha\}$ be orbit
representatives of the $d$-cells ($d \ge 0$) and write $S_\alpha$
for $\stab(e_\alpha) \in \F$. Define
$$
    C_d = \bigoplus_\alpha M \left( G/S_\alpha \right) \, .
$$
If $g e'$ is a typical $(d-1)$-cell in the boundary of $e_\alpha$
then $g^{-1} \stab(e_\alpha) g \subset \stab(e')$, giving a
$G$-map (write $S'$ for $\stab(e')$) $$\phi: G/S_\alpha
\rightarrow G/S'\, ,$$ which induces a homomorphism $M(\phi)
\colon M \left( G/S_\alpha\right) \rightarrow M\left(G/S'\right)$.
This yields a differential $\partial_d \colon C_d \rightarrow
C_{d-1}$, and the Bredon homology groups $H^\F_i\left(X;M\right)$
correspond to the homology of $(C_*, \partial_*)$. Bredon
cohomology is defined analogously, for $M$ a right Bredon module.

\subsection{Proper actions and the representation ring}
We are interested in the case $X = \EG$, $\F = \Fin(G)$ ---the
family of all finite subgroups of $G$--- and $M = \mathcal{R}$ the
complex representation ring, considered as a Bredon module as
follows. On objects we set
    $$
        \mathcal{R}(G/K) = R_\mathbb{C}(K), \quad K \in \Fin(G)
    $$
the ring of complex representations of the finite group $K$
(viewed just as an abelian group). For a $G$-map $\phi: G/K
\rightarrow G/L$ we have $g^{-1}Kg \subset L$ for some $g \in G$
so define $\mathcal{R}(\phi): R_\mathbb{C}(K) \rightarrow
R_\mathbb{C}(L)$ as induction from $g^{-1}Kg$ into $L$ once we
identify $R_\mathbb{C}(g^{-1}Kg)$ with $R_\mathbb{C}(K)$.

We state two useful results about the zero degree and the higher
degree Bredon homology groups $H_i^\Fin\left( \EG; \RR \right)$,
and a K\"unneth formula.
\begin{proposition}\label{prop:RankH0}
Let $G$ be group and denote by \textup{FC($G$)} the set of conjugacy
classes of elements of finite order in $G$. Then there is an
isomorphism
$$
H_0^\mathfrak{Fin} \left(\EG; \mathcal{R}\right)
\otimes_\Z \mathbb{C} \cong \mathbb{C}[\textup{FC}(G)]\, .
$$
\end{proposition}
\begin{proof} See Definition 3.18 and Theorem 3.19 in Mislin's
notes \cite{MislinValette:ProperGroupActions}.
\end{proof}
\noindent Consequently, the rank of the zeroth Bredon homology
group above coincides with the number of conjugacy classes of
elements of finite order in $G$.

Define the \emph{singular set} $X^{sing}$ of a $G$-set $X$ as the
subspace of the points with non-trivial stabilizers. The following
is Lemma 3.21 in Mislin's notes
\cite{MislinValette:ProperGroupActions}.
\begin{proposition}\label{prop:EGsing}
Let $G$ be an arbitrary group. Then
there is a natural map
$$ H_i^{\mathfrak{Fin}}\left( \EG{}; \mathcal{R} \right)
\longrightarrow H_i(\underline{B}G; \mathbb{Z})\, , $$ which is an
isomorphism in dimensions $i > \textrm{dim}
\left(\EG{}^{\textit{sing}}\right) + 1$ and injective in dimension
$i = \textrm{dim} \left(\EG{}^{\textit{sing}}\right) + 1$.
\end{proposition}

There is a K\"unneth formula for the direct product of two groups.
Given a group $G$, define $H^\Fin_i\left(G;\RR\right)$ as
$H^{\mathfrak{Fin}(G)}_i(X;\RR^G)$, where $X$ is any model of
\EG{} and $\RR^G$ is the representation ring as a Bredon module
over $\mathcal{O}_{\mathfrak{Fin}(G)}G$
(cf.~\cite{MislinValette:ProperGroupActions}).

\begin{proposition}\label{thm:KunnethFormulaGH}
Let $G$ and $H$ be two groups. For every $n \geq 0$ there is a
split exact sequence
        \begin{eqnarray*}
        0 \rightarrow \bigoplus_{i+j=n} \Big( H^{\mathfrak{Fin}}_i (G;\RR) \otimes
        H^{\mathfrak{Fin}}_j (H;\RR) \Big) \rightarrow
        H_n^{\mathfrak{Fin}}
        (G \ttimes H; \RR) \rightarrow \\
        \bigoplus_{i+j=n-1} \textup{Tor} \Big(
        H^{\mathfrak{Fin}}_i
        (G;\RR),H^{\mathfrak{Fin}}_j (H;\RR)\Big) \rightarrow 0\, .
        \end{eqnarray*}
\end{proposition}
\begin{proof} Details can be found in
\cite[\S 3]{SanchezGarcia:BredonHomCoxGroups} (see also
\cite{GG99}).
\end{proof}

\subsection{Equivariant $K$-homology}\label{section:EquivariantKhomology}
There is an equivariant version of $K$-homology, denoted
$K^G_i(-)$ and defined in \cite{DL98} (see also \cite{Joa03})
using spaces and spectra over the orbit category of $G$. It was
originally defined in \cite{BCH94} using Kasparov's $KK$-theory.
We will only recall the properties we need.

Equivariant $K$-homology satisfies Bott mod-2 periodicity, so we
only consider $K^G_0$ and $K^G_1$. For any subgroup $H$ of $G$ we
have
$$
    K^G_i\left(G/H\right) = K_i\left(C^*_r(H)\right) \, ,
$$
that is, its value at one-orbit spaces corresponds to the
$K$-theory of the reduced $C^*$-algebra of the typical stabilizer.
If $H$ is a finite subgroup then $C^*_r(H) = \mathbb{C}H$ and
$$
    K_i^G\left(G/H\right) = K_i(\mathbb{C}H) = \left\{
    \begin{array}{cl} R_\mathbb{C}(H) & i = 0\, ,\\ 0 & i = 1 \, .
    \end{array} \right.
$$
This allows us to view $K_i^G(-)$ as a Bredon module over
$\mathcal{O}_\Fin G$.

We can use an equivariant Atiyah-Hirzebruch spectral sequence to
compute the $K^G$-homology of a proper $G$-CW-complex $X$ from its
Bredon homology (see
\cite[pp.~49-50]{MislinValette:ProperGroupActions} for details),
as
$$
    E^2_{p,q} = H^\Fin_p\left( X; K^G_q(-) \right) \Rightarrow
    K^G_{p+q}\left(X\right)\, .
$$
In the simple case when Bredon homology concentrates at low degree
we deduce the following fact.
\begin{proposition} \label{prop:LowDegreeBredonHomology}Write $H_i = H^\Fin_i(X;\mathcal{R})$ and
$K_i = K^G_i(X)$. If $H_i = 0$ for $i \ge 2$ then $K_0 = H_0$ and
$K_1 = H_1$.
\end{proposition}

\section{A model for \underline{\textup{E}}\,\SL}
\subsection{The symmetric space}
We describe a first model of the classifying space
$\underline{E}\SLn$, for any $n \ge 2$
(cf.~\cite[pp.~38-40]{Brown:Cohomology}). Let $Q(n)$ be the space
of real, symmetric, positive definite $n\ttimes n$ matrices
(equivalently, positive definite quadratic forms on
$\mathbb{R}^n$). Multiplication by positive scalars gives an
action whose quotient $X(n) = Q(n)/\mathbb{R}^+$ is called the
\emph{symmetric space}. Since the right action of $GL(n,\R)$ on
$Q(n)$ given by $A\cdot g = g^t A g$ is transitive, with typical
stabilizer $O(n)$, we can identify $Q(n)$ with $GL(n,\R)/O(n)$. On
the other hand, the inclusion of $SL(n,\mathbb{R})/SO(n)$ into
$GL(n,\R)/O(n)$ is an \SL{}-equivariant homotopy equivalence, also
up to multiplication by positive scalars on the codomain.
Therefore, the symmetric space $X(n)$ is \SL{}-equivariant
homotopy equivalent to $SL(n,\mathbb{R})/SO(n)$.

The action of $GL(n,\mathbb{R})$ on $Q(n)$ induces and action on
$X(n)$, which restricts to a $\SLn$-action. As a $\SLn$-space,
$X(n)$ has finite stabilizers, since \SLn{} is a discrete subgroup
of $GL(n,\R)$ and $O(n)$ is compact. Moreover, we have the
following.
    \begin{proposition} The symmetric space $X(n) \simeq SL(n,\R)/SO(n)$
    is a model of $\underline{E}\SLn$, of dimension $n(n+1)/2 -
    1$.
    \end{proposition}
    \begin{proof}
The space $Q(n)$ is clearly a (convex) cone, so contractible.
Since the action is linear, the fixed point subspace $Q(n)^H$, $H
\le \SLn$, is also a cone. It is not empty whenever $H$ is a
finite subgroup; take for instance the `average point' of the
orbit of any $A \in Q(n)$,
$$
    \frac{1}{|H|} \sum_{h \in H} A \cdot h \; \in Q(n)^H\, .
$$
Note that $A \in Q(n)^H$ if and only if its class $[A] \in
X(n)^H$, so that $X(n)^H$ is neither empty if $H$ if finite.
Moreover, if we fix a representative of each class (for instance,
choose a matrix norm $\|\cdot\|$ and define $A_u = A/\|A\|$), then
there is a contracting homotopy $H([B], t)$ $=$ $[tA_u + (1-t)
B_u]$ for any fixed $[A] \in X^H$.
    \end{proof}
\begin{remark}
More generally, if $\Gamma$ is a discrete subgroup of a Lie group
$G$ with finitely many connected components, take $K$ a maximal
compact subgroup of $G$. Then $G/K$ is a model for \EG{} and,
therefore, a model for $\underline{E}\Gamma$ \cite{BCH94}.
\end{remark}

\subsection{Deformation retractions}
There is a better model of $\underline{E}\SLn$, obtained as an
$\SLn$-equivariant deformation retract of the symmetric space
$X(n)$, of dimension $n(n-1)/2$. It is obtained via reduction
theory of quadratic forms (see
\cite[pp.~213-17]{Brown:Cohomology}). However, finding an explicit
cellular decomposition and describing the stabilizers is in
general quite laborious. For \SL{} this has been done by Soul\'e
in \cite{Sou78} and also by Henn in \cite{Henn99}. We now describe
(without proofs) the deformation retract and orbit space for
$n=3$, following Soul\'e.

Denote the elements of $Q = Q(3)$, respectively of $X = X(3)$, as
$$
    A = (a_{ij}) \in Q \, ,\quad [A] = \{\lambda A \ | \ \lambda \in
    \mathbb{R}^+\} \in X \, .
$$
Recall the right action of $\SL$ on $Q$, $A\cdot g = g^tAg$, which
extends to an action on $X$. Each orbit in $X$ has a total
preorder as follows. Given $A = B \cdot g$ we say that $[A] < [B]$
if the sequence of diagonal elements of $A$ is smaller than the
one of $B$ with respect to the lexicographic order in
$\mathbb{R}^3$. (This is well-defined : if $\lambda A = B \cdot
g'$ then $\lambda B = B \cdot (g' g^{-1})$ so $\lambda = 1$ and $g
= g'$.) An element $[A] \in X$ is called \emph{reduced} if it is
minimal in its orbit. Define the subspace
$$
    Y = \{ [A] \textrm{ reduced and } a_{11} = a_{22} = a_{33}
    \}\, .
$$
    \begin{proposition} The space $Y$ is an $\SL$-deformation retract of $X$
    and therefore a model of $\underline{E}\SL$, of dimension 3.
    \end{proposition}
    \begin{proof}
    The result follows from Theorem 1 in \cite{Sou78}.
    \end{proof}
\begin{remark}
The minimal dimension for a model of $\underline{E}\SL$ is
actually three. In general, one can prove that the strict upper
triangular group in \SLn{} has \emph{cohomological dimension}
$n(n-1)/2$, so any model of $\underline{E}\SLn$ has dimension at
least that; see \cite[Chp.~VIII]{Brown:Cohomology}.
\end{remark}

\subsection{Description of the orbit space}
In this section we give an equivariant cellular decomposition and
describe the stabilizers for the orbit space $\SL\backslash Y =
\underline{B}\SL$. We follow the approach and notation in
Soul\'e's paper \cite{Sou78}, although an equivalent and detailed
description of this model can be found in Henn's \cite{Henn99}.

\begin{figure}[hbt]
    \begin{center}
    \includegraphics[scale = 0.6]{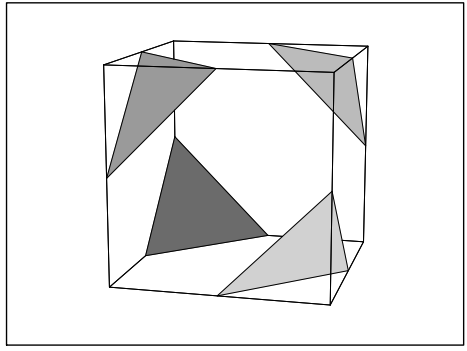}
    \caption{Truncated cube $C$}
    \label{fig:TruncatedCube}
    \end{center}
\end{figure}

\begin{figure}[hbt]
    \begin{center}
    \includegraphics[scale = 0.6]{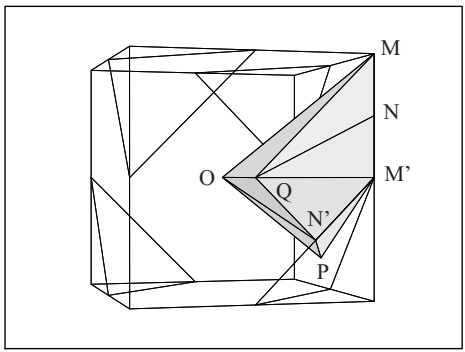}
    \caption{Triangulation of the fundamental domain $E$}
    \label{fig:SouleFundamentalDomain}
    \end{center}
\end{figure}

Let $C$ be the truncated cube of $\R^3$ with centre $(0,0,0)$ and
side length 2, truncated at the vertices $(1,1,-1)$, $(1,-1,1)$,
$(-1,1,1,)$ and $(-1,-1,-1)$, through the mid-points of the
corresponding sides (Figure \ref{fig:TruncatedCube}). %
An element in $Y$ can be uniquely written as the class of a matrix
$$
    A = \left(%
    \begin{array}{ccc}
      2 & z & y \\
      z & 2 & x \\
      y & x & 2 \\
    \end{array}%
    \right) = A(x,y,z)\, ,
$$
so that we can identify $[A]$ with the point $(x,y,z) \in
\mathbb{R}^3$. It can be shown that $A(x,y,z)$ is reduced if and
only if $(x,y,z)\in C$. Let $E$ be the subspace of $C$ given by
the points $(x,y,z)$ satisfying
\begin{equation*} \label{eqn:DefiningEqunsForFundDomainE}
    \begin{array}{c}
    |z| \leq y \leq x \leq 1\\
    z - x - y + 2 \ge 0 \, .
    \end{array}
\end{equation*}
We can give an explicit triangulation of $E$ as shown in Figure
\ref{fig:SouleFundamentalDomain}. The vertices are
    $$
    \begin{array}{ccccccc}
    O  & = & (0,0,0) & \hspace{5pt} & Q & = & (1,0,0)\\
    M  & = &  (1,1,1) & \hspace{5pt} & N  & = &  (1,1,1/2) \\
    M'  & = &  (1,1,0) & \hspace{5pt} & N'  & = &  (1,1/2,-1/2)\\
    P   & = &  (2/3,2/3,-2/3)\, . & & & & \\
    \end{array}
    $$
Note that the elements of \SL
    $$
    q_1 = \begin{pmatrix}
        1 & 0 & 0 \\
        0 & 1 & 1 \\
        0 & 0 & -1 \\
      \end{pmatrix} \textrm{ and }\,
    q_2 = \begin{pmatrix}
        -1 & 0 & 0 \\
        0 & 1 & 1 \\
        0 & 0 & -1 \\
      \end{pmatrix}
    $$
send $M, N, Q$ to $M', N', Q$ and $N, N', M', Q$ to $N', N, M', Q$
respectively. Consequently, we must identify the following
simplices in the quotient
    $$ \begin{array}{l}
    M \equiv M' \, ,\  N \equiv N' \, ,\ QM \equiv QM' \, ,\ QN \equiv QN' \, ,\ MN \equiv
    M'N \equiv M'N'\\ \textrm{and } \ QMN \equiv QM'N \equiv QM'N' \,. \end{array}
    $$
This identifications correspond to folding over the triangles
$QMN$, $QNM'$ and $QM'N'$ along the edges $QN$ and $QM'$
respectively.

\begin{theorem}
The space $E$ with the identifications above is a 3-dimensional
model of $\textrm{\underline{B}}\,\SL$.
\end{theorem}
\begin{proof}
See Theorem 2 in \cite{Sou78} or Theorem 2.4 in \cite{Henn99}.
\end{proof}

We now describe the orbits of cells and corresponding stabilizers.
This can be found in Theorem 2 of Soul\'e's article (although we
use a cellular structure instead of a simplicial one) or in
Section 2.5 of Henn's work. We have changed the chosen generators
so that they agree with the presentations on page
\pageref{tbl:CharTables}. We summarize the information in the
following table.

\vspace{12pt}

\label{table:InfoSL(3,Z)}
\begin{tabular}{lllc|lllc} %
    \hline
    \multicolumn{4}{c|}{Vertices}               &    \multicolumn{4}{c}{2-cells}\\
    \hline
    $v_1$ & $O$  & $g_2$, $g_3$  & $S_4$        &    $t_1$ & $OQM$   & $g_2$    & $C_2$   \\
    $v_2$ & $Q$  & $g_4$, $g_5$  & $D_6$        &    $t_2$ & $QM'N$  & $g_1$    & $\{1\}$     \\
    $v_3$ & $M$  & $g_6$, $g_7$  & $S_4$        &    $t_3$ & $M'N'P$ & $g_{12}$, $g_{14}$   & $C_2 \ttimes C_2$\\
    $v_4$ & $N$  & $g_6$, $g_8$  & $D_4$        &    $t_4$ & $OQN'P$ & $g_5$  & $C_2$   \\
    $v_5$ & $P$  & $g_5$, $g_9$ & $S_4$         &    $t_5$ & $OMM'P$ & $g_6$    & $C_2$   \\
    %\cline{1-4} \multicolumn{4}{c|}{Edges}     &    \\
    %\cline{1-4}
    \hline
    \multicolumn{4}{c|}{Edges}                  &    \multicolumn{4}{|c}{3-cell}\\
    \hline
    $e_1$ & $OQ$    & $g_2$, $g_5$      & $C_2 \ttimes C_2$  &  $T_1$ &  & $g_1$ & $\{1\}$\\
    $e_2$ & $OM$    & $g_6$, $g_{10}$     & $D_3$  \\
    $e_3$ & $OP$    & $g_6$, $g_5$      & $D_3$  \\
    $e_4$ & $QM$    & $g_2$             & $C_2$  \\
    $e_5$ & $QN'$  &  $g_5$             & $C_2$  \\
    $e_6$ & $MN$    & $g_6$, $g_{11}$   & $C_2 \ttimes C_2$\\
    $e_7$ & $M'P$   & $g_6$, $g_{12}$      & $D_4$\\
    $e_8$ & $N'P$   & $g_5$, $g_{13}$   & $D_4$ \\

    \end{tabular}

    \vspace{12pt}

The first column is an enumeration of equivalence classes of
cells; the second lists a representative of each class; the third
column gives generating elements for the stabilizer of the given
representative; and the last one is the isomorphism type of the
stabilizer. We use the following notation: $\{1\}$ denotes the
trivial group, $C_n$ the cyclic group of $n$ elements, $D_n$ the
dihedral group with $2n$ elements and $S_n$ the symetric group of
permutations on $n$ objects. The generating elements referred to
above are
$$
\begin{array}{ccc}
g_1 = \begin{pmatrix}  1&0&0 \\ 0&1&0 \\ 0&0&1  \end{pmatrix} &

g_2 = \begin{pmatrix}  -1&0&0 \\ 0&0&-1 \\ 0&-1&0  \end{pmatrix} &
g_3 = \begin{pmatrix}  0&0&1 \\ 0&1&0 \\ -1&0&0  \end{pmatrix}
\\[2em]
g_4 = \begin{pmatrix}  -1&0&0 \\ 0&1&1 \\ 0&0&-1  \end{pmatrix} &
g_5 = \begin{pmatrix}  -1&0&0 \\ 0&0&1 \\ 0&1&0  \end{pmatrix} &
g_6 = \begin{pmatrix}  0&-1&0 \\ -1&0&0 \\ 0&0&-1  \end{pmatrix}
\\[2em]
g_7 = \begin{pmatrix}  0&0&-1 \\ -1&0&0 \\ 1&1&1  \end{pmatrix} &
g_8 = \begin{pmatrix}  -1&0&0 \\ 0&1&0 \\ 0&-1&-1  \end{pmatrix} &
g_9 = \begin{pmatrix}  0&0&-1 \\ -1&0&-1 \\ 0&1&1  \end{pmatrix}
\\[2em]
g_{10} = \begin{pmatrix}  0&0&-1 \\ 0&-1&0 \\ -1&0&0
\end{pmatrix} & g_{11} = \begin{pmatrix}  -1&0&0 \\ 0&-1&0 \\
1&1&1\end{pmatrix} & g_{12} = \begin{pmatrix}  0&-1&-1 \\ 0&-1&0
\\ -1&1&0
\end{pmatrix}\\[2em]
g_{13} = \begin{pmatrix}  0&1&1 \\ 1&0&1 \\ 0&0&-1  \end{pmatrix}
& g_{14} =\begin{pmatrix}  -1&0&0 \\ -1&0&-1 \\ 1&-1&0
\end{pmatrix}\, .
    \end{array}
    $$
Finally, the cell boundaries can be easily read from Figure
\ref{fig:SouleFundamentalDomain}, once we fix an orientation;
namely, the ordering of the vertices $O < Q < M < M' < N < N' < P$
induces an orientation in $E$ and also in $\underline{B}\SL$ =
$E/\!\sim$. Thus, the boundaries with respect to the orbit
representatives are (recall the identifications given by $q_1$ and
$q_2$)
$$
\begin{array}{llll}
\partial e_1 = v_2 - v_1 & \partial e_2 = v_3 - v_1 &
    \partial e_3 = v_5 - v_1 &
\partial e_4 = v_3 - v_2 \\ \partial e_5 = v_4\cdot q_2 - v_2 &
    \partial e_6 = v_4 - v_3 &
\partial e_7 = v_5 - v_3\cdot q_1 & \partial e_8 = v_5 - v_4\cdot q_2\\[0.8em]

\multicolumn{2}{l}{
  \partial t_1 =  e_1 - e_2 + e_4} & \multicolumn{2}{l}{
                \partial t_2 =  e_4\cdot q_1 - e_5\cdot q_2 + e_6\cdot q_1q_2}\\
\multicolumn{2}{l}{  \partial t_3 =   e_6\cdot q_1 - e_7 + e_8} &
\multicolumn{2}{l}{  \partial t_4 = e_1 - e_3 + e_5 + e_8}\\
\multicolumn{2}{l}{  \partial t_5 =   e_2 - e_3 + e_6 - e_6\cdot
q_1q_2 + e_7} & \multicolumn{2}{l}{
                \partial T_1 = -t_1 + t_2 - t_3 + t_4 - t_5\, .}\\
\end{array}
$$
We now have all the information needed to compute the
corresponding Bredon homology, which is carried out in the next
section.

\section{Bredon homology of $\underline{E}\SL$}

The Bredon chain complex associated to $Y = \underline{E}\SL$ is
\begin{eqnarray}\label{eqn:BredonChainComplexESL(3,Z)}
    \lefteqn{0 \longrightarrow R_{\mathbb{C}}\left(\stab(T_1)\right)
    \stackrel{\Psi_3}{\longrightarrow}
    \bigoplus_{i=1}^5 R_{\mathbb{C}}\left(\stab(t_i)\right) \stackrel{\Psi_2}{\longrightarrow}}\nonumber\\
    & & \bigoplus_{j=1}^8 R_{\mathbb{C}}\left(\stab(e_j)\right) \stackrel{\Psi_1}{\longrightarrow}
    \bigoplus_{k=1}^5 R_{\mathbb{C}}\left(\stab(v_k)\right)
    \longrightarrow
    0 \, ,
\end{eqnarray}
where the $\Psi_i$'s are given by induction among representation
rings: if the boundary of a $d$-cell $e^d$ is, in terms of orbit
representatives of $(d-1)$-cells,
$$ \partial e^d = \sum_{j=1}^n e^{d-1}_j\cdot g_j $$
and $\tau \in R_{\mathbb{C}}(\stab(e^d))$, then $\Psi_d(\tau) =
\tau \uparrow \stab(e^{d-1}_j)\, ,$ where $\uparrow$ represents
induction with respect to $g\cdot\stab(e^d)\cdot g^{-1} \subset
\stab(e^{d-1}_j)$. Note that we write $\rho\uparrow H$ for
induction into a supergroup $H$ and $\rho\downarrow H$ for
restriction into a subgroup (we omit the group when it is clear
from the context).

To compute $\Psi_d(\tau)$, we use two basic facts:
\begin{itemize}
\item any representation (or character) can be uniquely written as
a sum of irreducible ones $ \tau = n_1 \rho_1 + \ldots + n_s
\rho_s$, with $n_i = (\tau | \rho_i)$ and $(\cdot | \cdot)$ the
usual scalar product of characters; \item Frobenius reciprocity:
$(\tau \uparrow | \rho)_H = (\tau | \rho\downarrow )_{H'}$ where
$H' \le H$.
\end{itemize}

Firstly, we note down the character tables of the groups appearing
as cell stabilizers. We write $\langle\, g \, \rangle$ for the
group generated by an element $g$, and $\langle \, gens \, |\,
rels \, \rangle$ for
a presentation of a group.\\%\\[1em]

\noindent\fbox{Trivial group}\quad $\{1\} = \langle \, g_1 \,
\rangle$ \nopagebreak \label{CharacterTables}
\begin{flushleft}
\begin{tabular}{c|c}
       $\{1\}$ & $g_1$ \\
       \hline
       $\tau$ & $1$\\
\end{tabular}
\end{flushleft}

\noindent\fbox{Cyclic group} \quad $C_2 = \langle\, g_i \,
\rangle$\nopagebreak
\begin{flushleft}
\begin{tabular}{c|rr}
       $C_2$ & $1$ & $g_i$ \\
       \hline
       $\rho_1$ & $1$ & $1$ \\
       $\rho_2$ & $1$ & $-1$ \\
\end{tabular}
\end{flushleft}

\label{tbl:CharTables}\noindent\fbox{Dihedral group} \quad $D_n =
\langle\, g_i, g_j \, \rangle = \langle \,g_i, g_j\,|\, (g_i)^2 =
(g_j)^2 = (g_ig_j)^n = 1\, \rangle$ \nopagebreak
\begin{flushleft}
\begin{tabular}{c|cc}
       $D_n$ & $(g_ig_j)^k$ & $g_j(g_ig_j)^k$ \\
       \hline
       $\chi_1$ & $1$ & $1$ \\
       $\chi_2$ & $1$ & $-1$ \\
       $\widehat{\chi_3}$ & $(-1)^k$ & $(-1)^k$ \\
       $\widehat{\chi_4}$ & $(-1)^k$ & $(-1)^{k+1}$ \\
       $\phi_p$ & $2\cos\left( 2\pi p k/n\right)$ & $0$ \\
\end{tabular}
\end{flushleft}
where $0 \le k \le n-1$, $p$ varies from 1 to $n/2-1$ ($n$ even)
or $(n-1)/2$ ($n$ odd) and the hat \ $\widehat{ }$ \ denotes a
character which only appears when $n$ is even. Note that $C_2
\times C_2 \cong D_2$.\\[1em]
\noindent\fbox{Symmetric group} \quad $S_4 = \langle\, g_i, g_j \,
\rangle$ where $g_i$ is a transposition and $g_j$ a cycle of
length 4. The character table, in cycle type notation, is
\begin{flushleft}
\begin{tabular}{c|rrrrr}
       $S_4$ & 1 & (12) & (123) & (1234) & (12)(34) \\
       \hline
       $\pi_1$ & $1$ & $1$ & 1 & 1 & 1 \\
       $\pi_2$ & $1$ & $-1$ & 1 & $-1$ & 1 \\
       $\pi_3$ & $2$ & $0$ & $-1$ & 0 & 2 \\
       $\pi_4$ & $3$ & $1$ & 0 & $-1$ & $-1$ \\
       $\pi_5$ & $3$ & $-1$ & 0 & 1 & $-1$ \\
\end{tabular}
\end{flushleft}
In the next sections we will need to agree on an ordering of the
conjugacy classes for each cell stabilizer. The same applies for
the irreducible characters. We fix both orderings as shown in the
character tables above. Note that for a dihedral group the
arranging of conjugacy classes (in terms of representatives) is
$$ \begin{array}{ll} n \textrm{ odd:} & 1,\, g_ig_j,\, (g_ig_j)^2,\,
\ldots,\,
(g_ig_j)^{(n-1)/2},\, g_j\\
n \textrm{ even:} & 1,\, g_ig_j,\, (g_ig_j)^2,\, \ldots,\,
(g_ig_j)^{n/2},\, g_j\,, g_j(g_ig_j). \end{array}
$$

Considering the ranks of the corresponding representation rings,
that is, the number of irreducible characters, we can view the
Bredon chain complex (\ref{eqn:BredonChainComplexESL(3,Z)}) as
(write $n \cdot \Z$ for $\Z^{\oplus n}$)
\begin{equation}\label{eqn:BredonChainComplexSL(3,Z)}
    \xymatrix{
    0 \ar[r] & \Z \ar[r]^-{\Psi_3} &
    11\cdot \Z \ar[r]^-{\Psi_2} &
    28\cdot \Z \ar[r]^-{\Psi_1} &
    26\cdot \Z \ar[r] & 0 \, .}
\end{equation}

\subsection{Computation of $\Psi_3$}
Denote by $\partial e$ the boundary of a $d$-cell $e$ in terms of
$(d-1)$-cells. Let $\tau$ be the trivial representation of
$\stab(T_1) = \{1\}$. We have
$$
    \partial T_1 = \sum_{i=1}^5
    (-1)^i t_i
\quad \Rightarrow \quad
    \Psi_3(\tau) = \sum_{i=1}^5 (-1)^i\, \tau \uparrow \stab(t_i)
    \, .
$$
Inducing the trivial representation gives the regular
representation $\rho_1 + \ldots + \rho_s$ so $\Psi_3$ is, written
as a matrix of an homomorphism of free abelian groups,

$$
    \Psi_3 (1) = \left( \begin{array}{cc|c|cccc|cc|cc}
    -1 & -1 &  1 &  -1 & -1 & -1 & -1 & 1 & 1 & -1 & -1\\
    \end{array}\right)\, .
$$
Here each group of $\pm1$s indicates the elements corresponding to
each representation ring $R_{\mathbb{C}}(\stab(t_i))$. This matrix
reduces by elementary operations to
$$ \Psi_3 \equiv
    \begin{pmatrix}
    1 & 0 & \cdots & 0
    \end{pmatrix}\, .
$$
In particular, $H_3 = \ker \Psi_3 = 0$ and $\textup{im}\, \Psi_3
\cong \Z$.

\subsection{Computation of $\Psi_2$}
For each 2-cell, we work out the induction map for each inclusion
(possibly after conjugation) of stabilizers. For $t_1$, we have
$\partial t_1 =  e_1 - e_2 + e_4$ so $\stab(t_1) \subset
\stab(e_i)$ for $i = 1, 2, 4$. The first inclusion is
$$
    C_2 = \stab(t_1) = \langle g_2 \rangle \subset \langle g_2,
    g_5 \rangle = \stab(e_1) = C_2 \times C_2 = D_2
    \, .
$$
Consider the irreducible characters $\rho_1$, $\rho_2$ in
$R_{\mathbb{C}}(\stab(t_1))$ and $\chi_1, \ldots, \chi_4$ in
$R_{\mathbb{C}}(\stab(e_1))$ as on the character tables on page
\pageref{CharacterTables}. Then
\begin{flushleft}
\begin{tabular}{c|rr|cc}
            & $1$ & $g_2$ & $(\rho_1\,|\,\chi_j\downarrow)$ & $(\rho_2\,|\,\chi_j\downarrow)$\\
       \hline
       $\chi_1\downarrow$ & $1$ & $1$ & 1 & 0\\
       $\chi_2\downarrow$ & $1$ & $-1$ & 0 & 1\\
       $\chi_3\downarrow$ & $1$ & $-1$ & 0 & 1\\
       $\chi_4\downarrow$ & $1$ & $1$ & 1 & 0
\end{tabular}
\end{flushleft}
Therefore,
$$
\begin{array}{rcl}
  \rho_1\uparrow &=& \chi_1 + \chi_4 \\
  \rho_2\uparrow &=& \chi_2 + \chi_3\, .
\end{array}
$$

\noindent The corresponding submatrix representing
$R_{\mathbb{C}}(\stab(t_1))
\stackrel{\textup{Ind}}{\longrightarrow} R_\mathbb{C}(\stab(e_1))$
is

$$
\begin{array}{cccccc}
% & & \multicolumn{4}{c}{\stab(e_1)}\\
  & & \chi_1 & \chi_2 & \chi_3 & \chi_4\\[0.3em]
%\multirow{2}{*}{$\stab(t_1)$}
 & \rho_1 \uparrow & 1 & 0 & 0 & 1\\
  & \rho_2 \uparrow & 0 & 1 & 1 & 0
\end{array}
$$

For the inclusion of $\stab(t_1) = \langle g_2 \rangle$ into
$\stab(e_2) = \langle g_6, g_{10} \rangle$, we have $g_2 =
g_{10}g_6g_{10}$, and
\begin{flushleft}
\begin{tabular}{c|rr|cc}
            & $1$ & $g_{10}(g_6g_{10}$ & $(\rho_1\,|\,\chi_j\downarrow)$ & $(\rho_2\,|\,\chi_j\downarrow)$\\
       \hline
       $\chi_1\downarrow$ & $1$ & $1$ & 1 & 0\\
       $\chi_2\downarrow$ & $1$ & $-1$ & 0 & 1\\
       $\phi_1\downarrow$ & $2$ & $0$ & 1 & 1\\
\end{tabular}
\end{flushleft}
Consequently,
$$
\begin{array}{rcl}
  \rho_1\uparrow &=& \chi_1 + \phi_1 \\
  \rho_2\uparrow &=& \chi_3 + \phi_1\, .
\end{array}
$$

The corresponding submatrix is (note the minus sign, since the
cell appears as $-e_2$ in $\partial t_1$)
$$
    \begin{array}{crrr}
     & \multicolumn{3}{c}{e_2}\\[0.25em]
     \multirow{2}{*}{$t_1$} & -1 & 0 & -1\\
                            & 0 & -1 & -1\\
     \end{array}
$$
Here we have simplified the notation; a matrix with top label $e$
and left label $t$ will represent the coefficients of the induced
characters of $R_\mathbb{C}\left(\stab(t)\right)$ into
$R_\mathbb{C}\left(\stab(e)\right)$, possibly after conjugation.

Now, $\stab(t_1) = \stab(e_4)$, so induction gives the identity
map. In brief:
$$
    \begin{array}{crrrr|rrr|rr}
     & \multicolumn{4}{c}{e_1} & \multicolumn{3}{c}{e_2} &
     \multicolumn{2}{c}{e_4}\\[0.25em]
     \multirow{2}{*}{$t_1$} & 1 & 0 & 0 & 1 & -1 & 0 & -1 & 1 & 0\\
                            & 0 & 1 & 1 & 0 & 0 & -1 & -1 & 0 & 1\\
     \end{array}
$$

\noindent The process is analogous for the other 2-cells. All the
inclusions among stabilizers are immediate except the ones listed
below.
$$
\begin{array}{rclclcl}
    \stab(t_3) & \subset & \stab(e_6\cdot q_1), &  & \stab(e_7), &  &
    \stab(e_8)\\
    g_{12} & = & q_1^{-1}g_6q_1 & = & g_{12} & = & g_{13}(g_5g_{13})^2\\
    g_{14} & = & q_1^{-1}g_{11}q_1 & = & g_{12}(g_6g_{12})^2 & = &
    (g_5g_{13})^2\, .\\[0.5em]
    \stab(t_5) & \subset & \stab\left(e_6\cdot q_1q_2\right) & &
    & & \\
    g_6 & = & (q_1q_2)^{-1}g_6g_{11}(q_1q_2)\, .
\end{array}
$$

We now show the results, which can be easily verified by a careful
reader.
$$
    \begin{array}{c cc|cc|cccc}
     & \multicolumn{2}{c}{e_4} & \multicolumn{2}{c}{e_5} &
     \multicolumn{4}{c}{e_6}\\[0.25em]
     \multirow{1}{*}{$t_2$} & 1 & 1 & -1 & -1 & 1 & 1 & 1 & 1\\
     \end{array}
$$
$$
    \begin{array}{c rrrr|rrrrr|rrrrr}
     & \multicolumn{4}{c}{e_6} & \multicolumn{5}{c}{e_7} &
     \multicolumn{5}{c}{e_8}\\[0.25em]
     \multirow{4}{*}{$t_3$} & 1 & 0 & 0 & 0  & -1 & 0 & -1 & 0 & 0  & 1 & 0 & 1 & 0 & 0\\
                            & 0 & 1 & 0 & 0  & 0 & -1 & 0 & -1 & 0  & 0 & 0 & 0 & 0 & 1\\
                            & 0 & 0 & 1 & 0  & 0 & 0 & 0 & 0 & -1   & 0 & 1 & 0 & 1 & 0\\
                            & 0 & 0 & 0 & 1  & 0 & 0 & 0 & 0 & -1   & 0 & 0 & 0 & 0 & 1\\
     \end{array}
$$
$$
    \begin{array}{c rrrr|rrr|rr|rrrrr}
     & \multicolumn{4}{c}{e_1} & \multicolumn{3}{c}{e_3} &
     \multicolumn{2}{c}{e_5} & \multicolumn{5}{c}{e_8}\\[0.25em]
     \multirow{2}{*}{$t_4$} & 1 & 0 & 1 & 0 & -1 & 0 & -1 & 1 & 0 & 1 & 0 & 0 & 1 & 1\\
                            & 0 & 1 & 0 & 1 & 0 & -1 & -1 & 0 & 1 & 0 & 1 & 1 & 0 & 1\\
     \end{array}
$$
$$
    \begin{array}{c rrr|rrr|rrrr|rrrrr}
     & \multicolumn{3}{c}{e_2} & \multicolumn{3}{c}{e_3} &
     \multicolumn{4}{c}{e_6} & \multicolumn{5}{c}{e_7}\\[0.25em]
     \multirow{2}{*}{$t_5$} & 1 & 0 & 1  & -1 & 0 & -1  & 0 & -1 & 0 & 1  & 1 & 0 & 0 & 1 & 1\\
                            & 0 & 1 & 1  & 0 & -1 & -1  & 0 & 1 & 0 & -1  & 0 & 1 & 1 & 0 & 1\\
     \end{array}
$$

\vspace{0.8em}

The submatrices above amount to an $11 \ttimes 28$ matrix
representing $\Psi_2$. This matrix can be reduced to its normal
form consisting of the identity of size 10 and zeroes elsewhere,
$$\Psi_2 \equiv \left(%
\begin{array}{c|c}
  Id_{10} & 0 \\
  \hline
  0 & 0 \\
\end{array}%
\right)\, .
$$

\subsection{Computation of $\Psi_1$}
The computations for $\Psi_1$ are similar and straightforward. The
relevant inclusions among stabilizers are the following. We give a
conjugacy representative as $(\sim g_i)$ when necessary.
$$
\begin{array}{rclcl}
    \stab(e_1) & \subset & \stab(v_1), &  & \stab(v_2)\\
    g_{2} & = & g_2 & = & g_{5}(g_4g_{5})^3\\
    g_{5} & = & g_{3}g_2g_3^{-1}g_2g_3 \,(\sim g_2)&  &\\[0.5em]
    \stab(e_2) & \subset & \stab(v_1), &  & \stab(v_3)\\
    g_{6} & = & g_{3}g_2g_3^{-1} & = & \\
    g_{10} & = & g_2g_3g_2g_3^{-1}g_2 \,(\sim g_2) & = &  g_7^{-1}g_6g_7\\[0.5em]
    \stab(e_3) & \subset & \stab(v_1), &  & \stab(v_5)\\
    g_{6} & = & g_{3}g_2g_3^{-1} & = & g_9^{-1}g_5g_9 \\
    g_{5} & = & g_3g_2g_3^{-1}g_2g_3 \,(\sim g_2) & & \\[0.5em]
    \stab(e_4) & \subset & \stab(v_2), &  & \stab(v_3)\\
    g_{6} & = & g_5(g_4g_5)^3 & = &  g_6g_7^{2}g_6g_7^{-1}\, (\sim g_6)
\end{array}$$
$$
\begin{array}{rclcl}
    \stab(e_5) & \subset & \stab(v_4\cdot q_2) &  & \\
    g_5 & = & q_2^{-1}g_8q_2 & & \\[0.5em]
    \stab(e_6) & \subset & \stab(v_3), &  & \stab(v_4)\\
    g_{11} & = & g_7g_6g_7^{-1}g_6g_7 \,(\sim g_6) & = & (g_6g_8)^2 \\[0.5em]
    \stab(e_7) & \subset & \stab(v_3\cdot q_1), &  & \stab(v_5)\\
    g_{6} & = & q_1^{-1}(g_{6}g_7^2g_6)q_1\,(\sim g_6) & = & g_9^{-1}g_5g_9 \\
    g_{12} & = & q_1^{-1}g_{6}q_1  & = & g_9^2 \,(\sim g_5) \\[0.5em]
    \stab(e_8) & \subset & \stab(v_4\cdot q_2), &  & \stab(v_5)\\
    g_{5} & = & q_2^{-1}g_{8}q_2 &  & \\
    g_{13} & = & q_2^{-1}g_{6}q_2  & = & g_5g_9^2g_5 \,(\sim g_5)\, .
\end{array}
$$

The matrices representing induction among stabilizers are the
following.
$$
    \begin{array}{c rrrrr|rrrrrr}
     & \multicolumn{5}{c}{v_1} & \multicolumn{6}{c}{v_2}\\[0.25em]
     \multirow{4}{*}{$e_1$} & -1 & 0 & -1 & -1 & 0     & 1 & 0 & 0 & 0 & 0 & 1 \\
                            & 0 & -1 & -1 & 0 & -1     & 0 & 1 & 0 & 0 & 0 & 1 \\
                            & 0 & 0 & 0 & -1 & -1      & 0 & 0 & 1 & 0 & 1 & 0 \\
                            & 0 & 0 & 0 & -1 & -1      & 0 & 0 & 0 & 1 & 1 & 0 \\
     \end{array}
$$
$$
    \begin{array}{c rrrrr|rrrrr}
     & \multicolumn{5}{c}{v_1} & \multicolumn{5}{c}{v_3}\\[0.25em]
     \multirow{4}{*}{$e_2$} & -1 & 0 & 0 & -1 & 0     & 1 & 0 & 0 & 1 & 0 \\
                            & 0 & -1 & 0 & 0 & -1     & 0 & 1 & 0 & 0 & 1 \\
                            & 0 & 0 & -1 & -1 & -1    & 0 & 0 & 1 & 1 & 1 \\
     \end{array}
$$
$$
    \begin{array}{c rrrrr|rrrrr}
     & \multicolumn{5}{c}{v_1} & \multicolumn{5}{c}{v_5}\\[0.25em]
     \multirow{4}{*}{$e_3$} & -1 & 0 & 0 & -1 & 0     & 1 & 0 & 0 & 1 & 0 \\
                            & 0 & -1 & 0 & 0 & -1     & 0 & 1 & 0 & 0 & 1 \\
                            & 0 & 0 & -1 & -1 & -1      & 0 & 0 & 1 & 1 & 1\\
     \end{array}
$$
$$
    \begin{array}{c rrrrrr|rrrrr}
     & \multicolumn{6}{c}{v_2} & \multicolumn{5}{c}{v_3}\\[0.25em]
     \multirow{4}{*}{$e_4$} & -1 & 0 & 0 & -1 & -1 & -1    & 1 & 0 & 1 & 2& 1 \\
                            & 0 & -1 & -1 & 0 & -1 & -1    & 0 & 1 & 1 & 1 & 2\\
     \end{array}
$$
$$
    \begin{array}{c rrrrrr|rrrrr}
     & \multicolumn{6}{c}{v_2} & \multicolumn{5}{c}{v_4}\\[0.25em]
     \multirow{4}{*}{$e_5$} & -1 & 0 & -1 & 0 & -1 & -1     & 1 & 0 & 1 & 0 & 1 \\
                            & 0 & -1 & 0 & -1 & -1 & -1     & 0 & 1 & 0 & 1 & 1 \\
     \end{array}
$$
$$
    \begin{array}{c rrrrr|rrrrr}
     & \multicolumn{5}{c}{v_3} & \multicolumn{5}{c}{v_4}\\[0.25em]
     \multirow{4}{*}{$e_6$} & -1 & 0 & -1 & -1 & 0     & 1 & 0 & 0 & 1 & 0 \\
                            & 0 & -1 & -1 & 0 & -1     & 0 & 0 & 0 & 0 & 1 \\
                            & 0 & 0 & 0 & -1 & -1      & 0 & 1 & 1 & 0 & 0 \\
                            & 0 & 0 & 0 & -1 & -1      & 0 & 0 & 0 & 0 & 1 \\
     \end{array}
$$
$$
    \begin{array}{c rrrrr|rrrrr}
     & \multicolumn{5}{c}{v_3} & \multicolumn{5}{c}{v_5}\\[0.25em]
     \multirow{4}{*}{$e_7$} & -1 & 0 & -1 & 0 & 0     & 1 & 0 & 1 & 0 & 0\\
                            & 0 & 0 & 0 & 0 & -1     & 0 & 0 & 0 & 0 & 1 \\
                            & 0 & 0 & 0 & -1 & 0      & 0 & 1 & 1 & 0 & 0 \\
                            & 0 & -1 & -1 & 0 & 0      & 0 & 0 & 0 & 1 & 0 \\
                            & 0 & 0 & 0 & -1 & -1      & 0 & 0 & 0 & 1 & 1\\
     \end{array}
$$
$$
    \begin{array}{c rrrrr|rrrrr}
     & \multicolumn{5}{c}{v_4} & \multicolumn{5}{c}{v_5}\\[0.25em]
     \multirow{4}{*}{$e_8$} & -1 & 0 & 0 & 0 & 0     & 1 & 0 & 1 & 0 & 0 \\
                            & 0 & -1 & 0 & 0 & 0     & 0 & 0 & 0 & 0 & 1 \\
                            & 0 & 0 & 0 & -1 & 0      & 0 & 1 & 1 & 0 & 0 \\
                            & 0 & 0 & -1 & 0 & 0      & 0 & 0 & 0 & 1 & 0 \\
                            & 0 & 0 & 0 & 0 & -1      & 0 & 0 & 0 & 1 & 1 \\

     \end{array}
$$

Altogether they form a $28 \ttimes 26$ matrix whose normal form is
$$\left(%
\begin{array}{c|c}
  Id_{18} & 0 \\
  \hline
  0 & 0 \\
\end{array}%
\right)\, .
$$
\begin{remark}
The computations have been verified with the help of a computer.
In fact, the author has implemented a program in GAP \cite{GAP},
which computes the Bredon homology with coefficients in the
representation ring of any finite proper $G$-CW-complex, from the
cell stabilizers and boundaries. Details and the code for the
algorithm can be found in \cite[Appendix A]{SanchezGarciaPhD}.
\end{remark}

We can now determine the Bredon homology of \SL{} from the chain
maps above. Note that if we have a exact sequence of free abelian
groups
$$
   \cdots \longrightarrow n \cdot \mathbb{Z} \stackrel{f}{\longrightarrow} m \cdot \mathbb{Z}
        \stackrel{g}{\longrightarrow} k \cdot \mathbb{Z}
        \longrightarrow \cdots
$$
with $f$ and $g$ represented by matrices $A$ and $B$ for some
fixed basis, then the homology at $m \cdot \mathbb{Z}$ is
\begin{equation*} \label{eqn:Ker/Im}
        \textrm{ker\,} (g) \,/\, \textrm{im\,} (f)
        \,\cong\, \mathbb{Z}/d_1\Z \oplus \ldots \mathbb{Z}/d_s\Z
        \oplus
        (m-s-r) \cdot \mathbb{Z} \, ,
\end{equation*}
where $r = \textrm{rank}\,(B)$ and $d_1, \ldots, d_s$ are the
elementary divisors of $A$. In our case, we obtain that the Bredon
homology of $\underline{E}\SL$ is
\begin{eqnarray*}
H^\Fin_0\left( \underline{E}\Gamma; \RR \right) &\cong& \Z^{\oplus 8}\\
H^\Fin_i\left( \underline{E}\Gamma; \RR \right) &=& 0 \quad
\forall \,i \neq 0\, .
\end{eqnarray*}
\begin{remark}
These results agree with the Bredon homology expected at degree 0
and 3. The dimension of the singular part of our model of
$\underline{E}\Gamma$ is 2, and $\underline{B}\Gamma$ is
contractible, so we have $H_i = 0$ for all $i \ge 3$, by
Proposition \ref{prop:EGsing}. On the other hand,
$\textup{rank}\,(H_0) = 8$ must be the number of conjugacy classes
of elements of finite order of \SL{}, by Proposition
\ref{prop:RankH0}; the latter can be deduced, for instance, from
the list of all finite subgroups (up to conjugacy) of \SL{} in
\cite{Tahara71}. Also, 8 is the alternating sum of the ranks in
the Bredon chain complex (\ref{eqn:BredonChainComplexSL(3,Z)}).
\end{remark}

\subsection{Equivariant $K$-homology}
Since the Bredon homology concentrates at degree 0, it coincides
with the $K^G$-homology (Proposition
\ref{prop:LowDegreeBredonHomology}), that is,
\begin{eqnarray*}
K^G_0\left( \underline{E}\Gamma \right) &=& \Z^{\oplus 8}\\
K^G_1\left( \underline{E}\Gamma \right) &=& 0\, .
\end{eqnarray*}
This amounts to the topological side of the Baum-Connes conjecture
and injects into the analytical side, that is, the $K$-theory of
$C^*_r(\SL)$.

\subsection{Results for \GL}
We have the direct product decomposition $\GL = \SL \times C_2$.
We can therefore use the K\"unneth formula for Bredon homology
from Proposition \ref{thm:KunnethFormulaGH}. Since $C_2$ is
finite, a one-point space is a model for $\underline{E}C_2$ and
its Bredon homology is $R_{\mathbb{C}}(C_2) \cong \Z^{\oplus 2}$
at degree 0 and vanishes elsewhere. Consequently,
$$
\begin{array}{rcl}
H^\Fin_0\left( \underline{E}\,\GL; \RR \right) &=& H^\Fin_0\left(
\underline{E}\SL; \RR \right) \otimes
    H^\Fin_0\left(\underline{E}C_2; \RR \right) \cong \Z^{\oplus 16}\\
H^\Fin_i\left( \underline{E}\,\GL; \RR \right) &=& 0 \, ,\quad
\textrm{if } i \neq 0\, .
\end{array}
$$
As before, these groups coincide with the equivariant
$K$-homology:
\begin{eqnarray*}
K^G_0\left( \underline{E}\,\GL\right) &\cong& \Z^{\oplus 16}\\
K^G_1\left( \underline{E}\,\GL\right) &=& 0\, .
\end{eqnarray*}

%%% THE BIBLIOGRAPHY %%%%%%%%%%%%%%%%%%%%%%%%%%%%%%%%%%%%%%%%%%%%
\bibliographystyle{amsplain}
\bibliography{MyBibliography}

\providecommand{\bysame}{\leavevmode\hbox to3em{\hrulefill}\thinspace}
\providecommand{\MR}{\relax\ifhmode\unskip\space\fi MR }
\providecommand{\MRhref}[2]{%
  \href{http://www.ams.org/mathscinet-getitem?mr=#1}{#2}
}
\providecommand{\href}[2]{#2}
\begin{thebibliography}{10}

\bibitem{BCH94}
P.~Baum, A.~Connes, and N.~Higson, \emph{Classifying spaces for proper actions
  and {$K$-theory} of group {$C^*$-algebras}}, Contemporary Mathematics
  \textbf{167} (1994), 241--291.

\bibitem{Brown:Cohomology}
K.~S. Brown, \emph{Cohomology of {Groups}}, Springer Graduate Text in
  Mathematics, no.~87, Springer-Verlag, 1982.

\bibitem{DL98}
J.~F. Davis and W.~L{\"u}ck, \emph{Spaces over a {Category} and {Assembly}
  {Maps} in {Isomorphism} {Conjectures} in {$K$}- and {$L$}-{Theory}},
  $K$-Theory \textbf{15} (1998), 201--252.

\bibitem{GAP}
The GAP~Group, \emph{{GAP -- Groups, Algorithms, and Programming, Version
  4.4}}, 2004, \verb+(http://www.gap-system.org)+.

\bibitem{GG99}
M.~Golasi{\'n}ski and D.~L. Gon{\c{c}}alves, \emph{Generalized
  {E}ilenberg-{Z}ilber type theorem and its equivariant applications}, Bull.
  Sci. Math. \textbf{123} (1999), no.~4, 285--298.

\bibitem{Henn99}
H.-W. Henn, \emph{The cohomology of {${\rm SL}(3,{\bf Z}[1/2])$}}, $K$-Theory
  \textbf{16} (1999), no.~4, 299--359.

\bibitem{HK97}
N.~Higson and G.~Kasparov, \emph{Operator {$K$}-theory for groups which act
  properly and isometrically on {H}ilbert space}, Electron. Res. Announc. Amer.
  Math. Soc. \textbf{3} (1997), 131--142 (electronic).

\bibitem{HLS02}
N.~Higson, V.~Lafforgue, and G.~Skandalis, \emph{Counterexamples to the
  {Baum}--{Connes} {Conjecture}}, Geom. Funct. Anal. \textbf{12} (2002), no.~2,
  330--354.

\bibitem{Joa03}
M.~Joachim, \emph{{$K$}-homology of {$C\sp \ast$}-categories and symmetric
  spectra representing {$K$}-homology}, Math. Ann. \textbf{327} (2003), no.~4,
  641--670.

\bibitem{Laf98}
V.~Lafforgue, \emph{Une d\'emonstration de la conjecture de {B}aum-{C}onnes
  pour les groupes r\'eductifs sur un corps {$p$}-adique et pour certains
  groupes discrets poss\'edant la propri\'et\'e ({T})}, C. R. Acad. Sci. Paris
  S\'er. I Math. \textbf{327} (1998), no.~5, 439--444.

\bibitem{LN01}
I.~J. Leary and B.~E.~A. Nucinkis, \emph{Every {CW-complex} is a classifying
  space for proper bundles}, Topology \textbf{40} (2001), 539--550.

\bibitem{MislinValette:ProperGroupActions}
G.~Mislin and A.~Valette, \emph{Proper group actions and the {B}aum-{C}onnes
  conjecture}, Advanced Courses in Mathematics. CRM Barcelona, Birkh\"auser
  Verlag, 2003.

\bibitem{SanchezGarcia:BredonHomCoxGroups}
R.~S\'anchez-Garc\'\i{}a, \emph{Equivariant {$K$-homology} for some {Coxeter}
  groups}, arXiv:math.KT/0604402.

\bibitem{SanchezGarciaPhD}
\bysame, \emph{Equivariant {$K$-homology} of the classifiying space for proper
  actions}, Ph.D. thesis, University of Southampton, 2005.

\bibitem{Sou78}
C.~Soul{\'{e}}, \emph{The cohomology of {$\textrm{SL}_3(\mathbb{Z})$}},
  Topology \textbf{17} (1978), 1--22.

\bibitem{Tahara71}
K.-i. Tahara, \emph{On the finite subgroups of {${\rm GL}(3,\,\bf Z)$}}, Nagoya
  Math. J. \textbf{41} (1971), 169--209.

\end{thebibliography}

%%%% END OF DOCUMENT %%%%%%%%%%%%%%%%%%%%%%%%%%%%%%%%%%%%%%%%%%%%
\end{document}